# Matrix-analytic solution of system of integral equations in three tandem servers

ABDULLAH ALKAFF *

**Abstract**

A matrix-analytic method is proposed for solving a system of linear integral equations arises in three tandem servers. The approach is by modelling the cumulative distribution function (cdf) of the service time as a matrix exponential function. The method transforms the system of linear integral equations into a system of linear algebraic equations, hence can produce closed-form solutions. Properties of the solution are discussed.

## 1 Introduction

This study is an attempt to devise a method for providing a general analytical solution to the following system of linear integral equations derived in [5]:

$$F_{I_2}(t) = 1 - \int_{x=0}^{\infty} F_{S_2}(x) F_{R_3}(x)\, dF_{S_1}(x+t) \qquad \text{for } t \geq 0 \tag{1}$$

$$F_{R_3}(t) = 1 - \int_{x=0}^{\infty} F_{I_2}(x)\, dF_{S_3}(x+t) \qquad \text{for } t \geq 0 \tag{2}$$

* Institute of Technology Sepuluh Nopember (ITS), Surabaya, Indonesia. alkaff@ee.its.ac.id.

These equations represent relationships among several random variables in a system comprises of three servers arranged in series without buffer space between two consecutive servers. The functions $F_{S_j}(x), j = 1, 2, 3,$ in these equations are known and represent the cumulative distribution function (cdf) of service time at server $j$ as the inputs to the system. The functions $F_{I_2}(x)$ and $F_{R_3}(x)$ represent cdf of two of the system's outputs that must be determined by solving the system of integral equations given the input $F_{S_j}(x), j = 1, 2, 3$. Other outputs' cdfs can be found directly once $F_{I_2}(x)$ and $F_{R_3}(x)$ are determined.

Solutions to these equations are given in [5], [6] and [4] for some simple model of $F_{S_j}(x)$. Furthermore, their solutions are only for the mean value of the most important variable for engineers, i.e. throughput rate. This paper presents a method for solving (1) and (2) in the case $F_{S_j}(x)$ has a more general model that covers all previously used cdfs. The solutions are functions, not just values, that can be used to determine the cdfs of all other outputs, hence contain versatile information regarding the behavior of the tandem servers beyond throughput rate. For this purpose, it is assumed that the inputs have matrix-based distributions. Matrix-based distribution is very general since it can be used to represent many scalar-based distributions, including those used in previous publications. On the other hand, using matrix-based distribution will produce simple matrix equations that can be solved algorithmically to find the cdfs of all variables in the system.

In mathematical term, matrix-based distribution basically is a scalar function represented as a matrix exponential function. It was used to provide a simple solution to the Volterra integral equation as a generalization of Pollaczek-Khinchin integral equation for waiting time in single server queueing systems [7]. In this study, the system of integral equations that need to be solved contains two or three matrix exponential functions. Analysis of systems with two or more related matrix-based distribution requires the use of Kronecker, or tensor, matrix operations. Such algebra has been employed for solving an integral equation in complex geometry [3]. Thus, the novelty of this study is that it proposes the use of Kronecker operations to solve a system of integral equations involving more than two matrix exponential functions.

## 2 Matrix-based distributions and Kronecker operations

Matrix-based distribution is introduced to model randomness more realistically than scalar exponential distribution yet has compact form and properties that can simplify many integral and algebraic operations. There are two popular matrix-based distributions known as phase-type (PH) distribution proposed in [8] and matrix-exponential (ME) distribution as a generalization of PH distribution introduced in [2].

PH distribution is based on the distribution of the time to absorption, $X$, in a finite, irreducible, homogenous continuous-time Markov process with transient states $\{1, 2, \dots, n\}$, and an absorbing state $(n + 1)$. The cdf of $X$ is given by

$$F_X(x) = 1 - \boldsymbol{a}e^{\boldsymbol{A}x}\boldsymbol{u} \qquad \text{for } x \geq 0 \tag{3}$$

where

- $\boldsymbol{a} = \big(a(1), a(2), \dots, a(n)\big)$ is the $1 \times n$ initial state probability vector, hence $0 \leq a(i) \leq 1$ for $i = 1,2, \dots, n$,
- $\boldsymbol{u}$ is a $n \times 1$ unity vector (vector whose elements equal unity), and
- $\boldsymbol{A}$ is $n \times n$ transition matrix among the transient states whose $A_{ij}$ entries represent the nonnegative transition rates from a transient state $i$ to a transient state $j$, its diagonal entries $A_{ii}$ must be negative, and the sum of all elements at the same row must be nonpositive. By this definition, $\boldsymbol{A}$ is nonsingular and $e^{\boldsymbol{A}x} \to 0$ as $x \to \infty$.

This function has a jump or discontinuity at $x = 0$ of size $a(n+1) = 1 - \boldsymbol{au}$. Hence, its derivative is given by

$$f_X(x) = -\boldsymbol{aA}e^{\boldsymbol{A}x}\boldsymbol{u} + (1 - \boldsymbol{au})\delta(x) \qquad \text{for } x \geq 0 \tag{4}$$

where $\delta(x)$ is an impulse function represents discontinuity of $F_X(x)$ at $x = 0$. The distribution $F_X(\cdot)$ is said to be of PH with representation $(n, \boldsymbol{a}, \boldsymbol{A})$ and written as a PH distribution. Denotes $X \backsim PH(n, \boldsymbol{a}, \boldsymbol{A})$ as a shorthand notation for 'random variable $X$ has PH distribution with PH representation $(n, \boldsymbol{a}, \boldsymbol{A})$.'

PH distribution has been generalized into ME distribution in [11] by allowing $\boldsymbol{a}, \boldsymbol{A}, \boldsymbol{u}$ to be arbitrary but must satisfy the requirements for (1) to be a valid cdf, i.e.

$$0 \leq \boldsymbol{au} \leq 1, \boldsymbol{A} \text{ is nonsingular, } -\boldsymbol{aAu} \geq 0, \text{ and } e^{\boldsymbol{A}x} \to 0 \text{ as } x \to \infty$$

which are automatically satisfies in PH distribution from its definition.

This study will be based on PH distribution because of its popularity. Nevertheless, the approach is applicable to ME distribution as well. Thus, cdf of service time at each station will be modelled as (3).

The use of matrix-based distribution in the analysis of stochastic systems requires matrix algebra operations known as Kronecker (or tensor) operations. Kronecker operations comprise the Kronecker product and the Kronecker sum, popular mathematical tools in matrix analysis [1].

**Definition 2.1**. If $\boldsymbol{A} = [A_{ij}]$ and $\boldsymbol{B} = [B_{ij}]$ are matrices of dimensions $m_1 \times m_2$ and $n_1 \times n_2$, then their Kronecker product $\boldsymbol{A} \otimes \boldsymbol{B}$ is a matrix of dimensions $n_1 m_1 \times n_2 m_2$, given in block-partitioned form as

$$\boldsymbol{A} \otimes \boldsymbol{B} = \begin{bmatrix} A_{11}\boldsymbol{B} & A_{12}\boldsymbol{B} & \dots & A_{1n_2}\boldsymbol{B} \\ \vdots & \vdots & \ddots & \vdots \\ A_{n_1 1}\boldsymbol{B} & A_{n_1 2}\boldsymbol{B} & \dots & A_{n_1 n_2}\boldsymbol{B} \end{bmatrix} \tag{5}$$

**Definition 2.2**. If $\boldsymbol{A}$ and $\boldsymbol{B}$ are matrices of dimensions $m \times m$ and $n \times n$, their Kronecker sum $\boldsymbol{A} \oplus \boldsymbol{B}$ is a matrix of dimension $mn \times mn$, defined as

$$\boldsymbol{A} \oplus \boldsymbol{B} = \boldsymbol{A} \otimes \boldsymbol{I}_n + \boldsymbol{I}_m \otimes \boldsymbol{B} \tag{6}$$

where $\boldsymbol{I}_m$ and $\boldsymbol{I}_n$ are identity matrices of dimensions $m \times m$ and $n \times n$, respectively.

The properties of Kronecker operations which are going to be employed in subsequent sections of this study are ($g$ and $h$ are scalars),

$$g \otimes h = gh \tag{7}$$

$$g \otimes \boldsymbol{A} = g\boldsymbol{A} \tag{8}$$

$$\boldsymbol{AB} \otimes \boldsymbol{CD} = (\boldsymbol{A} \otimes \boldsymbol{C})(\boldsymbol{B} \otimes \boldsymbol{D}) \tag{9}$$

$$e^{\boldsymbol{A}t} \otimes e^{\boldsymbol{B}t} = e^{(\boldsymbol{A} \oplus \boldsymbol{B})t} \tag{10}$$

The above properties can be extended to more complex expressions involving several Kronecker products by using the distributive property of the Kronecker product. However, the Kronecker product and the Kronecker sum are not commutative.

**Result 2.1**. $\left(\boldsymbol{a}_1 e^{\boldsymbol{A}_1 t}\boldsymbol{u}_1\right)\left(\boldsymbol{a}_2 e^{\boldsymbol{A}_2 t}\boldsymbol{u}_2\right) = (\boldsymbol{a}_1 \otimes \boldsymbol{a}_2) e^{(\boldsymbol{A}_1 \oplus \boldsymbol{A}_2)t}(\boldsymbol{u}_1 \otimes \boldsymbol{u}_2)$

The proof follows (7),

$$\left(\boldsymbol{a}_1 e^{\boldsymbol{A}_1 t}\boldsymbol{u}_1\right)\left(\boldsymbol{a}_2 e^{\boldsymbol{A}_2 t}\boldsymbol{u}_2\right) = \left(\boldsymbol{a}_1 e^{\boldsymbol{A}_1 t}\boldsymbol{u}_1\right) \otimes \left(\boldsymbol{a}_2 e^{\boldsymbol{A}_2 t}\boldsymbol{u}_2\right)$$

and use (9) repeatedly

$$\left(\boldsymbol{a}_1 e^{\boldsymbol{A}_1 t}\boldsymbol{u}_1\right) \otimes \left(\boldsymbol{a}_2 e^{\boldsymbol{A}_2 t}\boldsymbol{u}_2\right) = (\boldsymbol{a}_1 \otimes \boldsymbol{a}_2)\left(e^{\boldsymbol{A}_1 t} \otimes e^{\boldsymbol{A}_2 t}\right)(\boldsymbol{u}_1 \otimes \boldsymbol{u}_2)$$

and then use (10) to combine the two exponential functions into one to complete the proof.

**Result 2.2**. $\boldsymbol{a} \otimes \boldsymbol{b} = \boldsymbol{a}(\boldsymbol{I}_m \otimes \boldsymbol{b})$ for $\boldsymbol{a}$ and $\boldsymbol{b}$ are $1 \times m$ and $1 \times n$ vectors, respectively.

The proof follows from the definition of the Kronecker product

$$\boldsymbol{a} \otimes \boldsymbol{b} = [a_1\boldsymbol{b}, a_2\boldsymbol{b}, \dots, a_m\boldsymbol{b}] = [a_1, a_2, \dots, a_m]\begin{bmatrix} \boldsymbol{b} & 0 & \dots & 0 \\ \vdots & \vdots & \ddots & \vdots \\ 0 & 0 & \dots & \boldsymbol{b} \end{bmatrix}$$

$$= \boldsymbol{a}\begin{bmatrix} 1 & \dots & 0 \\ \vdots & \ddots & \vdots \\ 0 & \dots & 1 \end{bmatrix} \otimes \boldsymbol{b} = \boldsymbol{a}(\boldsymbol{I}_m \otimes \boldsymbol{b})$$

**Result 2.3** [1]: If $\boldsymbol{A}_1$and $\boldsymbol{A}_2$ are nonsingular then $\boldsymbol{A}_1 \oplus \boldsymbol{A}_2$ is nonsingular.

## 3 The distribution of $\boldsymbol{R}_3$ and $\boldsymbol{I}_2$

In the following development, PH distributions will be used to model the cdfs of $S_1$ and $S_3$, i.e. $F_{S_1}(t)$ and $F_{S_3}(t)$, as inputs to the system.

Let $S_1 \backsim PH(n_1, \boldsymbol{a}_1, \boldsymbol{A}_1)$ and $S_3 \backsim PH(n_3, \boldsymbol{a}_3, \boldsymbol{A}_3)$,

$$F_{S_1}(t) = 1 - \boldsymbol{a}_1 e^{\boldsymbol{A}_1 t} \boldsymbol{u}_1 \quad t \geq 0 \tag{11}$$

$$F_{S_3}(t) = 1 - \boldsymbol{a}_3 e^{\boldsymbol{A}_3 t} \boldsymbol{u}_3 \quad t \geq 0 \tag{12}$$

Their derivatives are pdf given in (4) written as

$$dF_{S_1}(t) = \left(-\boldsymbol{a}_1' e^{\boldsymbol{A}_1 t} \boldsymbol{u}_1 + (1 - \boldsymbol{a}_1 \boldsymbol{u}_1)\delta(t)\right) dt \tag{13}$$

$$dF_{S_3}(t) = \left(-\boldsymbol{a}_3' e^{\boldsymbol{A}_3 t} \boldsymbol{u}_3 + (1 - \boldsymbol{a}_3 \boldsymbol{u}_3)\delta(t)\right) dt \tag{14}$$

where $\boldsymbol{a}_1' = \boldsymbol{a}_1 \boldsymbol{A}_1$ and $\boldsymbol{a}_3' = \boldsymbol{a}_3 \boldsymbol{A}_3$.

The cdf of $S_2$ is assumed to be arbitrary but is Laplace transformable.

Substituting (14) into (1) yields

$$F_{R_3}(t) = 1 - \int_{x=0}^{\infty} F_{I_2}(x) \left(-\boldsymbol{a}_3' e^{\boldsymbol{A}_3(x+t)} \boldsymbol{u}_3 + (1 - \boldsymbol{a}_3 \boldsymbol{u}_3)\delta(x+t)\right) dx$$

$$F_{R_3}(t) = 1 - \int_{x=0}^{\infty} F_{I_2}(x)(-\boldsymbol{a}_3' e^{\boldsymbol{A}_3 x} e^{\boldsymbol{A}_3 t}\, \boldsymbol{u}_3)\, dx - (1 - \boldsymbol{a}_3 \boldsymbol{u}_3) \int_{x=0}^{\infty} F_{I_2}(x)\delta(x+t)\, dx$$

$$F_{R_3}(t) = 1 - \left(\int_{x=0}^{\infty} F_{I_2}(x)(-\boldsymbol{a}_3' e^{\boldsymbol{A}_3 x} dx)\right) e^{\boldsymbol{A}_3 t} \boldsymbol{u}_3$$

The second integration is zero since the impulse occurs outside the range of integration. Hence, $F_{R_3}(t)$ can be written as

$$F_{R_3}(t) = 1 - \boldsymbol{r}_3 e^{\boldsymbol{A}_3 t} \boldsymbol{u}_3 \tag{15}$$

where $\boldsymbol{r}_3$ is a row vector of order $n_3$ given by

$$\boldsymbol{r}_3 = -\int_0^{\infty} F_{I_2}(x) \boldsymbol{a}_3'\, e^{\boldsymbol{A}_3 x} dx \tag{16}$$

**Result 3.1**. $R_3 \backsim PH(n_3, \boldsymbol{r}_3, \boldsymbol{A}_3)$

To prove this result is sufficient by showing that $\boldsymbol{r}_3$ satisfy the requirement for $F_{R_3}(t)$ to be a valid cdf. Since $F_{I_2}(x)$ is nonnegative, and $-\int_0^\infty \boldsymbol{a}_3'\, e^{\boldsymbol{A}_3 x} dx = -\int_0^\infty \boldsymbol{a}_3 \boldsymbol{A}_3\, e^{\boldsymbol{A}_3 x} dx = \boldsymbol{a}_3$, hence $\boldsymbol{r}_3$ is nonnegative. Multiplying (16) with $\boldsymbol{u}_3$, results in

$$\boldsymbol{r}_3\boldsymbol{u}_3 = -\int_0^\infty F_{I_2}(x)\boldsymbol{a}_3' e^{\boldsymbol{A}_3 x}\boldsymbol{u}_3\, dx < -\int_0^\infty \boldsymbol{a}_3' e^{\boldsymbol{A}_3 x}\boldsymbol{u}_3\, dx = \boldsymbol{a}_3\boldsymbol{u}_3 \le 1$$

Thus, each element of $\boldsymbol{r}_3$ is between 0 and 1 which satisfy the requirement of PH distribution.

Substitute (13) into (2) to yield

$$F_{I_2}(t) = 1 - \int_0^\infty F_{S_2}(x)F_{R_3}(x)\left(-\boldsymbol{a}_1' e^{\boldsymbol{A}_1(x+t)}\boldsymbol{u}_1 + \delta(x+t)\right) dx$$

$$F_{I_2}(t) = 1 - \left(\int_0^\infty F_{S_2}(x)F_{R_3}(x)\left(-\boldsymbol{a}_1' e^{\boldsymbol{A}_1 x}\right)dx\right) e^{\boldsymbol{A}_1 t}\boldsymbol{u}_1 - \left(\int_0^\infty F_{S_2}(x)F_{R_3}(x)\delta(x+t)dx\right)$$

The second integration is also zero for the same reason that the impulse occurs outside the range of integration. Hence,

$$F_{I_2}(t) = 1 - \left(\int_0^\infty F_{S_2}(x)F_{R_3}(x)\left(-\boldsymbol{a}_1' e^{\boldsymbol{A}_1 x}\right)dx\right) e^{\boldsymbol{A}_1 t}\boldsymbol{u}_1$$

which can be written as

$$F_{I_2}(t) = 1 - \boldsymbol{d}_2 e^{\boldsymbol{A}_1 t}\boldsymbol{u}_1 \tag{17}$$

where $\boldsymbol{d}_2$ is a row vector of order $n_1$, the same order as $\boldsymbol{A}_1$, given by

$$\boldsymbol{d}_2 = -\int_0^\infty F_{S_2}(x)F_{R_3}(x)\boldsymbol{a}_1' e^{\boldsymbol{A}_1 x} dx \tag{18}$$

**Result 3.2**. $I_2 \backsim PH(n_1, \boldsymbol{d}_2, \boldsymbol{A}_1)$

This claim can be proved in the same manner as proving Result 3.1.

**Result 3.3**. PH distribution is closed in the system governed by equations (1) and (2) irrespective of the distribution of $S_2$.

The proof immediately follows Results 3.1 and 3.2 as they represent the only two outputs of the system.

## 4 The system of linear equations for $\boldsymbol{r}_3$ and $\boldsymbol{d}_2$

One advantage of using matrix-based distributions is that it can transform integral equations arise in stochastic systems into algebraic equations. Equations (1) and (2) therefore can be transformed into a system of algebraic equations when matrix-based distributions are used as inputs. As will be clear later, such transformation does not require any specific assumption on the distribution of the service period of the middle server.

Substituting (17) into (16) yields,

$$\boldsymbol{r}_3 = -\int_0^\infty \left(1 - \boldsymbol{d}_2 e^{\boldsymbol{A}_1 x} \boldsymbol{u}_1\right)\left(\boldsymbol{a}_3' e^{\boldsymbol{A}_3 x}\right) dx$$

For convenience, express $\boldsymbol{a}_3' e^{\boldsymbol{A}_3 x}$ as $\boldsymbol{a}_3' e^{\boldsymbol{A}_3 x} \boldsymbol{I}_{n_3}$ where $\boldsymbol{I}_{n_3}$ is an $n_3 \times n_3$ identity matrix. Hence,

$$\boldsymbol{r}_3 = -\int_0^\infty \boldsymbol{a}_3' e^{\boldsymbol{A}_3 x} dx + \int_0^\infty \left(\boldsymbol{d}_2 e^{\boldsymbol{A}_1 x} \boldsymbol{u}_1\right)\left(\boldsymbol{a}_3' e^{\boldsymbol{A}_3 x} \boldsymbol{I}_{n_3}\right) dx$$

Use (6) to replace the multiplication of a scalar $\boldsymbol{d}_2 e^{\boldsymbol{A}_1 x} \boldsymbol{u}_1$ and a vector $\boldsymbol{a}_3' e^{\boldsymbol{A}_3 x} \boldsymbol{I}_{n_3}$ with the Kronecker product to obtain

$$\boldsymbol{r}_3 = -\int_0^\infty \boldsymbol{a}_3' e^{\boldsymbol{A}_3 x} dx + \int_0^\infty \left(\boldsymbol{d}_2 e^{\boldsymbol{A}_1 x} \boldsymbol{u}_1\right) \otimes \left(\boldsymbol{a}_3' e^{\boldsymbol{A}_3 x} \boldsymbol{I}_{n_3}\right) dx$$

Rearranging the Kronecker products using the property given in (10) yields

$$\boldsymbol{r}_3 = \boldsymbol{a}_3 + \int_0^\infty (\boldsymbol{d}_2 \otimes \boldsymbol{a}_3')(e^{\boldsymbol{A}_1 x} \otimes e^{\boldsymbol{A}_3 x})(\boldsymbol{u}_1 \otimes \boldsymbol{I}_{n_3}) dx$$

$$\boldsymbol{r}_3 = \boldsymbol{a}_3 + (\boldsymbol{d}_2 \otimes \boldsymbol{a}_3') \int_0^\infty e^{(\boldsymbol{A}_1 \oplus \boldsymbol{A}_3)x} dx \, \left(\boldsymbol{u}_1 \otimes \boldsymbol{I}_{n_3}\right)$$

$$\boldsymbol{r}_3 = \boldsymbol{a}_3 - (\boldsymbol{d}_2 \otimes \boldsymbol{a}_3')(\boldsymbol{A}_1 \oplus \boldsymbol{A}_3)^{-1}(\boldsymbol{u}_1 \otimes \boldsymbol{I}_{n_3})$$

(as $\int_0^\infty e^{(\boldsymbol{A}_1 \oplus \boldsymbol{A}_3)x} dx = -(\boldsymbol{A}_1 \oplus \boldsymbol{A}_3)^{-1}$ )

To simplify the notation, define $\boldsymbol{A}_{13} = \boldsymbol{A}_1 \oplus \boldsymbol{A}_3$ and $\boldsymbol{A}_{31} = \boldsymbol{A}_3 \oplus \boldsymbol{A}_1$. Hence

$$\boldsymbol{r}_3 = \boldsymbol{a}_3 - (\boldsymbol{d}_2 \otimes \boldsymbol{a}_3')(\boldsymbol{A}_{13})^{-1}(\boldsymbol{u}_1 \otimes \boldsymbol{I}_{n_3}) \qquad (19)$$

Use Result 2.2 to convert (19) into

$$\boldsymbol{r}_3 = \boldsymbol{a}_3 - \boldsymbol{d}_2\left(\boldsymbol{I}_{n_1} \otimes \boldsymbol{a}_3'\right)(\boldsymbol{A}_{13})^{-1}(\boldsymbol{u}_1 \otimes \boldsymbol{I}_{n_3}) \tag{20}$$

which can be simplified into

$$\boldsymbol{r}_3 = \boldsymbol{a}_3 - \boldsymbol{d}_2\boldsymbol{B} \tag{21}$$

where

$$\boldsymbol{B} = \left(\boldsymbol{I}_{n_1} \otimes \boldsymbol{a}_3'\right)(\boldsymbol{A}_{13})^{-1}(\boldsymbol{u}_1 \otimes \boldsymbol{I}_{n_3}) \tag{22}$$

Result 2.3 guarantees the nonsingularity of $\boldsymbol{A}_{13}$, hence the uniqueness of $\boldsymbol{B}$.

Similarly, substitute (15) into (18) to yield

$$\boldsymbol{d}_2 = -\int_0^\infty F_{S_2}(x)\,(1 - \boldsymbol{r}_3 e^{\boldsymbol{A}_3 x}\boldsymbol{u}_3)\left(\boldsymbol{a}_1' e^{\boldsymbol{A}_1 x}\right)dx$$

$$\boldsymbol{d}_2 = -\boldsymbol{a}_1'\int_0^\infty F_{S_2}(x)e^{\boldsymbol{A}_1 x}dx + \int_0^\infty F_{S_2}(x)(\boldsymbol{r}_3 e^{\boldsymbol{A}_3 x}\boldsymbol{u}_3)(\boldsymbol{a}_1' e^{\boldsymbol{A}_1 x}\boldsymbol{I}_{n_1})dx$$

where $\boldsymbol{I}_{n_1}$ is an $n_1 \times n_1$ identity matrix added for convenience

Employing the same technique for the last integral as that used for the derivation of (20) yields

$$\boldsymbol{d}_2 = -\boldsymbol{a}_1'\int_0^\infty F_{S_2}(x)e^{\boldsymbol{A}_1 x}dx + \int_0^\infty F_{S_2}(x)(\boldsymbol{r}_3 e^{\boldsymbol{A}_3 x}\boldsymbol{u}_3) \otimes (\boldsymbol{a}_1' e^{\boldsymbol{A}_1 x}\boldsymbol{I}_{n_1})\,dx$$

$$= -\boldsymbol{a}_1' F_{S_2}^*(-\boldsymbol{A}_1) + \int_0^\infty F_{S_2}(x)(\boldsymbol{r}_3 \otimes \boldsymbol{a}_1')\left(e^{\boldsymbol{A}_3 x} \otimes e^{\boldsymbol{A}_1 x}\right)(\boldsymbol{u}_3 \otimes \boldsymbol{I}_{n_1})dx$$

$$= -\boldsymbol{a}_1' F_{S_2}^*(-\boldsymbol{A}_1) + (\boldsymbol{r}_3 \otimes \boldsymbol{a}_1')(\int_0^\infty F_{S_2}(x)e^{(\boldsymbol{A}_3 \oplus \boldsymbol{A}_1)\,x}dx)(\boldsymbol{u}_3 \otimes \boldsymbol{I}_{n_1})$$

$$\boldsymbol{d}_2 = -\boldsymbol{a}_1' F_{S_2}^*(-\boldsymbol{A}_1) + (\boldsymbol{r}_3 \otimes \boldsymbol{a}_1')F_{S_2}^*(-\boldsymbol{A}_{31})(\boldsymbol{u}_3 \otimes \boldsymbol{I}_{n_1}) \tag{23}$$

where $F_{S_2}^*(\mathrm{s})$ is the Laplace transform of $F_{S_2}(t)$ defined as $\int_0^\infty F_{S_2}(t)e^{-st}dt$.

Use Result 2.2 to convert (23) into

$$\boldsymbol{d}_2 = -\boldsymbol{a}_1' F_{S_2}^*(-\boldsymbol{A}_1) + \boldsymbol{r}_3\left(\boldsymbol{I}_{n_3} \otimes \boldsymbol{a}_1'\right)F_{S_2}^*(-\boldsymbol{A}_{31})(\boldsymbol{u}_3 \otimes \boldsymbol{I}_{n_1}) \tag{24}$$

which can be simplified into

$$\boldsymbol{d}_2 = -\boldsymbol{a}_1' F_{S_2}^*(-\boldsymbol{A}_1) + \boldsymbol{r}_3\boldsymbol{C} \tag{25}$$

where

$$\boldsymbol{C} = \left(\boldsymbol{I}_{n_3} \otimes \boldsymbol{a}_1'\right) F_{S_2}^*(-\boldsymbol{A}_{31})(\boldsymbol{u}_3 \otimes \boldsymbol{I}_{n_1}) \quad (26)$$

The existence of $F_{S_2}^*(-\boldsymbol{A}_{31})$ follows from the assumption that $F_{S_2}(t)$ is Laplace transformable. The case $S_2$ has PH distribution will be discussed in Section 4.5. Consequently, the uniqueness of $\boldsymbol{C}$ is guaranteed.

Equations (21) and (25) constitute a system of linear algebraic equations for $\boldsymbol{r}_3$ and $\boldsymbol{d}_2$,

$$\boldsymbol{d}_2\boldsymbol{B} + \boldsymbol{r}_3\boldsymbol{I}_{n_3} = \boldsymbol{a}_3 \quad (27)$$

$$\boldsymbol{d}_2\boldsymbol{I}_{n_1} - \boldsymbol{r}_3\boldsymbol{C} = -\boldsymbol{a}_1' F_{S_2}^*(-\boldsymbol{A}_1) \quad (28)$$

which can be combined as

$$[\boldsymbol{d}_2, \boldsymbol{r}_3]\begin{bmatrix} \boldsymbol{B} & \boldsymbol{I}_{n_1} \\ \boldsymbol{I}_{n_3} & -\boldsymbol{C} \end{bmatrix} = [\boldsymbol{a}_3, -\boldsymbol{a}_1' F_{S_2}^*(-\boldsymbol{A}_1)] \quad (29)$$

**Result 4.1**. PH distribution transforms the system of integral equations given in (1) and (2) into a system of linear algebraic equations given in (27) and (28).

It is instructive to note that the solution to the system of integral equations can be derived without making any specific assumption on the distribution of the middle server. It holds for any $F_{S_2}(t)$ that is Laplace transformable. Furthermore, $F_{R_3}(t)$ and $F_{I_2}(t)$, as the solutions of the system of integral equations, are both PH distribution irrespective of the distribution underlying $F_{S_2}(t)$.

## 5 The closed-form solutions

Finding $\boldsymbol{r}_3$ and $\boldsymbol{d}_2$ can be achieved by solving (29). It can be simplified further by directly substituting (25) into (21) to obtain

$$\boldsymbol{r}_3 = \boldsymbol{a}_3 - \left(-\boldsymbol{a}_1' F_{S_2}^*(-\boldsymbol{A}_1) + \boldsymbol{r}_3\boldsymbol{C}\right)\boldsymbol{B}$$

or

$$\boldsymbol{r}_3(\boldsymbol{I}_{n_3} + \boldsymbol{C}\boldsymbol{B}) = \boldsymbol{a}_3 + \boldsymbol{a}_1' F_{S_2}^*(-\boldsymbol{A}_1)\boldsymbol{B} \quad (30)$$

Equation (30) can be written simply as

$$\boldsymbol{r}_3\boldsymbol{G} = \boldsymbol{g} \tag{31}$$

or

$$\boldsymbol{r}_3 = \boldsymbol{g}\boldsymbol{G}^{-1} \tag{32}$$

where

$$\boldsymbol{G} = \boldsymbol{I}_{n_3} + \boldsymbol{C}\boldsymbol{B} \tag{33}$$

$$\boldsymbol{g} = \boldsymbol{a}_3 + \boldsymbol{a}_1' F_{S_2}^*(-\boldsymbol{A}_1)\boldsymbol{B} \tag{34}$$

Equation (32) represents a closed-form formula for $\boldsymbol{r}_3$. The closed-form formula for $F_{R_3}(t)$ is readily available since its PH representation is $(n_3, \boldsymbol{r}_3, \boldsymbol{A}_3)$. The closed-form formula for $\boldsymbol{d}_2$ can be obtained in the same fashion by eliminating $\boldsymbol{r}_3$ from (21) and (25).

The number of linear equations in (32) is equal to the order of PH representation of the service time of the third server, i.e. $n_3$. Another way of simplifying (29) is by eliminating $\boldsymbol{r}_3$ from (21) and (25) to yield a system of linear equations for $\boldsymbol{d}_2$. The order of such equations is equal to the order of PH representation of the first server, i.e. $n_1$. If the complexity of the problem is measured by the number of linear equations to be solved, then it is equal to the minimum of the orders of the first and last servers.

**Result 5.1**. The complexity of the problem is min $(n_1, n_3)$

## 6 Laplace transform of matrix functions

Evaluating $F_{S_2}^*(-\boldsymbol{A}_1)$ is straightforward when $F_{S_2}(t)$ represents scalar-based cdf. However, when it is given in a matrix-based cdf, its Laplace transform must be evaluated differently. Let $S_2 \backsim PH(n_2, \boldsymbol{a}_2, \boldsymbol{A}_2)$, or

$$F_{S_2}(t) = 1 - \boldsymbol{a}_2 e^{\boldsymbol{A}_2 t}\boldsymbol{u}_2 \tag{35}$$

$$F_{S_2}^*(-\boldsymbol{A}_1) = \int_0^\infty e^{\boldsymbol{A}_1 x} F_{S_2}(x)\, dx = \int_0^\infty e^{\boldsymbol{A}_1 x}\left(1 - \boldsymbol{a}_2 e^{\boldsymbol{A}_2 x}\boldsymbol{u}_2\right) dx$$

$$= \int_0^\infty e^{A_1 x}\, dx - \int_0^\infty e^{A_1 x}(\boldsymbol{a}_2 e^{A_2 x}\boldsymbol{u}_2)\, dx$$

To combine the two exponential functions of $x$ in the last integral together, the multiplication of matrix $e^{A_1 x}$ with scalar $\boldsymbol{a}_2 e^{A_2 x}\boldsymbol{u}$ is represented as the Kronecker product. For simplicity, the term $e^{A_1 x}$ is first multiplied by the identity matrix of the same dimension as $\boldsymbol{A}_1$. Hence,

$$F^*_{S_2}(-\boldsymbol{A}_1) = -\boldsymbol{A}_1^{-1} - \int_0^\infty \left(\boldsymbol{I}_{n_1} e^{A_1 x} \boldsymbol{I}_{n_1}\right) \otimes \left(\boldsymbol{a}_2 e^{A_2 x}\boldsymbol{u}_2\right) dx$$

$$F^*_{S_2}(-\boldsymbol{A}_1) = -\boldsymbol{A}_1^{-1} - \left(\boldsymbol{I}_{n_1} \otimes \boldsymbol{a}_2\right)\left(\int_0^\infty \left(e^{(A_1 \oplus A_2)\, x}\right) dx\right)\left(\boldsymbol{I}_{n_1} \otimes \boldsymbol{u}_2\right)$$

$$F^*_{S_2}(-\boldsymbol{A}_1) = -\boldsymbol{A}_1^{-1} + \left(\boldsymbol{I}_{n_1} \otimes \boldsymbol{a}_2\right)(\boldsymbol{A}_{12})^{-1}\left(\boldsymbol{I}_{n_1} \otimes \boldsymbol{u}_2\right) \tag{36}$$

where $\boldsymbol{A}_{12} = \boldsymbol{A}_1 \oplus \boldsymbol{A}_2$.

Result 2.3 guarantees the nonsingularity of $\boldsymbol{A}_{12}$, and PH distribution guarantees the nonsingularity of $\boldsymbol{A}_1$. Hence, both guarantee the existence and uniqueness of $F^*_{S_2}(-\boldsymbol{A}_1)$.

## 7 Numerical results

As an example, assume that service times at all servers have PH distributions with the following representations:

– PH representation of $S_1$:

$$\boldsymbol{a}_1 = [0.5 \quad 0.3]\ ;\ \ \boldsymbol{A}_1 = \begin{bmatrix} -1.0330 & 0.3099 \\ 0.3984 & -1.3281 \end{bmatrix}$$

– PH representation of $S_2$:

$$\boldsymbol{a}_2 = [0.5, 0.3, 0.2]\ ;\ \ \boldsymbol{A}_2 = \begin{bmatrix} -1.6321 & 0.8161 & 0.8161 \\ 0 & -32643 & 0 \\ 2.9379 & 2.4482 & -4.8964 \end{bmatrix}$$

– PH representation of $S_3$:

$$\boldsymbol{a}_1 = [1, 0, 0, 0]\ ;\ \ \boldsymbol{A}_1 = \begin{bmatrix} -4 & 4 & 0 & 0 \\ 0 & -4 & 4 & 0 \\ 0 & 0 & -4 & 4 \\ 0 & 0 & 0 & -4 \end{bmatrix}$$

Note that cdf of service time of $S_1$ has a jump of size $1 - 0.5 - 0.3 = 0.2$ at $t = 0$.

Matrix $\boldsymbol{B}$ and $\boldsymbol{C}$ are calculated using (22) and (26) to produce

$$\boldsymbol{B} = \begin{bmatrix} 0.8449 & -0.1325 & -0.1127 & -0.0956 \\ 0.8139 & -0.1496 & -0.1207 & -0.0978 \end{bmatrix}; \quad \boldsymbol{C} = \begin{bmatrix} -0.0853 & -0.0514 \\ -0.0598 & -0.0361 \\ -0.0351 & -0.0212 \\ -0.0137 & -0.0083 \end{bmatrix}$$

Matrix $\boldsymbol{G}$ and vector $\boldsymbol{g}$ are calculated using (33) and (34) to yield

$$\boldsymbol{G} = \begin{bmatrix} 0.8862 & 0.0190 & 0.0158 & 0.0132 \\ -0.0799 & 1.0133 & 0.0111 & 0.0093 \\ -0.0469 & 0.0078 & 1.0065 & 0.0054 \\ -0.0184 & 0.0031 & 0.0026 & 1.0021 \end{bmatrix}; \quad \boldsymbol{g} = [0.6315, 0.0614, 0.0512, 0.0426]$$

The values of $\boldsymbol{r}_3$ and $\boldsymbol{d}_2$ calculated by solving (29) to yield

$$\boldsymbol{r}_3 = [0.7196 \quad 0.0467 \quad 0.0389 \quad 0.0324]; \quad \boldsymbol{d}_2 = [0.2111 \quad 0.1253]$$

The value of $\boldsymbol{r}_3$ can also be obtained using (32) which will produce the same value. Using their PH representations, cdfs of these random variables are plotted as presented in Figure 1.

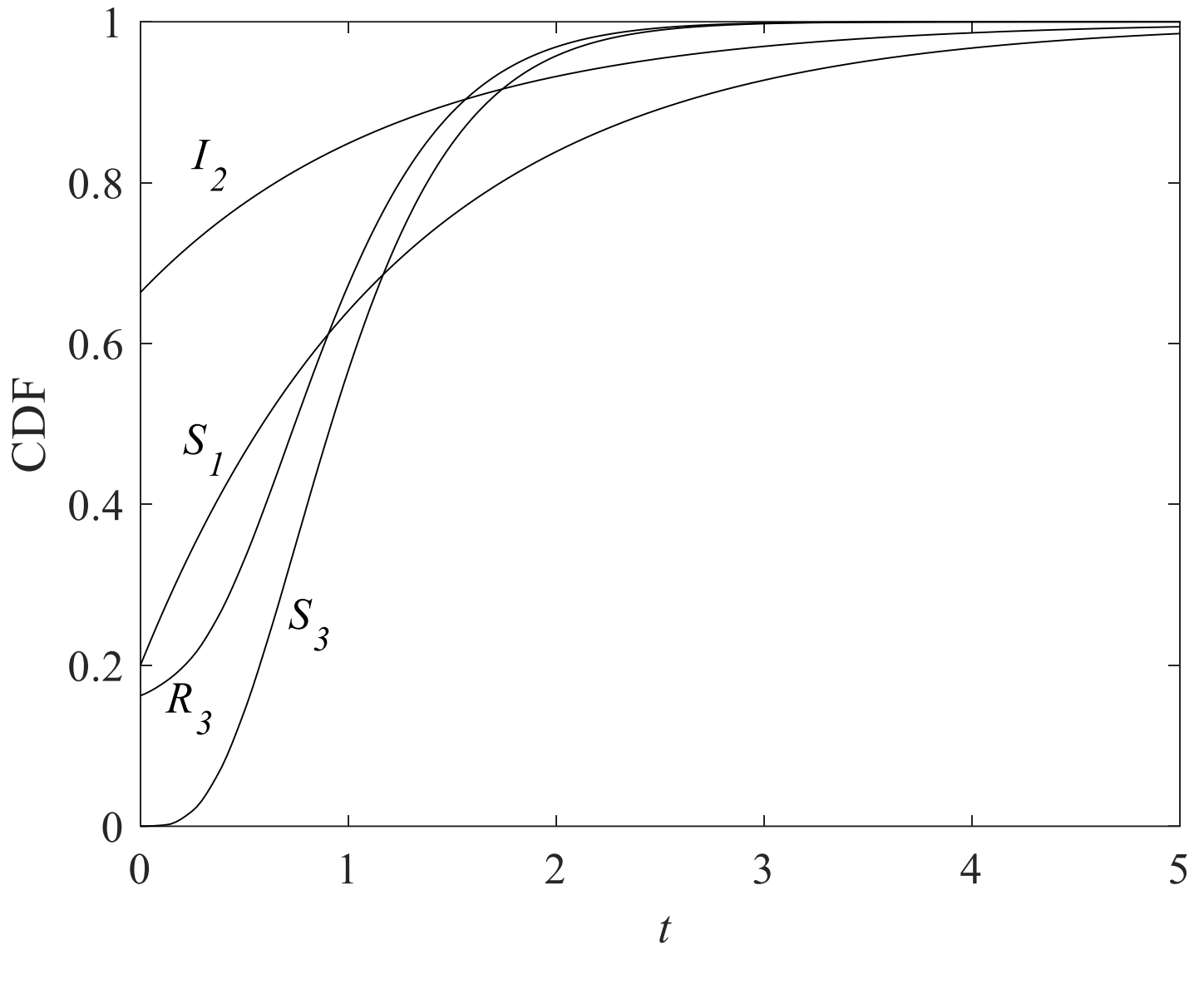

Figure 1: The graphs of cdfs of $S_1, S_3, R_3$, and $I_2$

## 8 Conclusions

Matrix-based distribution combined with Kronecker operations have been shown useful for solving a system of linear integral equations arises in three tandem servers. The solutions are closed-form formulas for cdfs of random variables in the equations.

An important result from this development is that PH distribution possesses closure properties in the system of linear integral equations, i.e. if the inputs to the system of linear integral equations have PH distributions then the outputs also have PH distributions. This result will open up opportunity, as future research, to derive closed-form formulas for all random variables in three tandem servers that are not present in the system of integral equations.

## References

[1] R. Bellman, *Introduction to Matrix Analysis*, (2nd Ed) Society for Industrial and Applied Mathematics, Philadelphia, 1997.

[2] M. Bladt and B. F. Nielsen, *Matrix-Exponential Distributions in Applied Probability*. Springer, 2017.

[3] E. Corona, A. Rahimian, and D. Zorin, "A Tensor-Train accelerated solver for integral equations in complex geometries," *J. Comput. Phys.*, 334 (2017), 145–169.

[4] S. Kekre, U. S. Rao, J. M. Swaminathan, and J. Zhang, "Reconfiguring a Remanufacturing Line at Visteon, Mexico," *Interfaces*, 33 (2003), 30–43.

[5] E. J. Muth, "Stochastic processes and their network representations associated with a production line queuing model," *Eur. J. Oper. Res.*, 15 (1984), 63–83.

[6] E. J. Muth and A. Alkaff, "The throughput rate of three-station production lines: A unifying solution," *Int. J. Prod. Res.*, 25 (1987),1405–1413

[7] M. F. Neuts, "Generalizations of the Pollaczek-Khinchin integral equation in the theory of queues," *Adv. Appl. Probab.*, 18 (1986), 952–990.

[8] M. F. Neuts, *Matrix-Geometric Solutions in Stochastic Models: An Algorithmic Approach*. Dover Publication, Inc., 1994.